\documentclass[11pt]{article}
\usepackage{amssymb}
\usepackage{mathrsfs}
\usepackage{amsmath}
\usepackage{cases}
\usepackage[dvips]{graphicx}
\author{Jun-Ming Zhu
 \footnote{E-mail address: jm\_zh@sohu.com, junming\_zhu@163.com}
}
\title{The solutions of
four $q$-functional equations
}
\date{}
\newtheorem{lem}{\quad\textbf{\Large Lemma}}[section]
\newtheorem{thm}[lem]{\quad\textbf{\Large Theorem}}
\begin{document}
 \maketitle
 \setcounter{section}{0}
 \textbf{Abstract} In this note we obtain the solutions of
four $q$-functional equations and express the solutions in
$q$-operator forms. These equations give sufficient conditions for
$q$-operator methods

\textbf{Key words:} Basic  hypergeometric series, $q$-Series,
$q$-Functional equation, $q$-Difference equation, q-Exponential
operator
\section{ Introduction }
~~~We follow the notation and terminology in \cite{gasper} , and for
convenience, we always assume that $0<|q|<1$. The $q$-shifted
factorials are defined by
\begin{equation*}
{(a;q)_n=}\left\{\begin{array}{ll}
1, &n=0, \\
(1-a)(1-aq)\cdots(1-aq^{n-1}), &n=1,2,3, \cdots\cdots,
\end{array}
\right.
\end{equation*}
The $q$-derivative operator is defined by
$$
D_q\{f(a)\}=\frac{f(a)-f(aq)}{a},
~~D_q^n\{f(a)\}=D_q\{D_q^{n-1}\{f(a)\}\}.
$$
The operator $\theta$ is defined by
$$
\theta=\eta^{-1}D_q, ~~\theta^n\{f(a)\}=\theta\{\theta^{n-1}\{f(a)\}\},
$$
where $\eta^{-1}\{f(a)\}=f(q^{-1}a).$ Both $D_q$ and $\theta$ are
obviously linear transforms, and by convention, $D_q^0$ and
$\theta^0$ are both understood as the identity.

The two $q$-exponential operators (see, \cite{chen-liu2,
chen-liu3,liu7,liu8,rogers}) are defined by
$$
T(bD_q)=\sum_{n=0}^{\infty}\frac{(bD_q)^{n}}{(q;q)_n},
$$
and
$$
 E(b\theta)=\sum_{n=0}^{\infty}\frac{q^{n(n-1)/
2}(b\theta)^{n}}{(q;q)_n},
$$
respectively.

The following operator was first introduced by Fang \cite{fang}. But
we follow Chen and Gu's notation \cite{chen-gu}. It seems to be more
convenient.
$$
T(a, b; D_q)=\sum_{n=0}^{\infty}\frac{(a, q)_n}{(q;q)_n}(bD_q)^{n}.
$$
Chen and Gu \cite{chen-gu} named $T(a,b;D_{q})$ as Cauchy operator.
But compared with the following operator, we called $T(a,b;D_{q})$
the first Cauchy operator in this paper.

The operator $T(a,b;\theta)$, introduced by Fang \cite{fang}, is
defined by
$$
T(a, b;\theta)=\sum_{n=0}^{\infty}\frac{(a,
q)_n}{(q;q)_n}(b\theta)^{n}.
$$
We called $T(a,b;\theta)$ the second Cauchy operator.

All the papers using the $q$-operator methods imply the following
definition but do not state it explicitly. Unless otherwise stated,
all the operators are applied with respect to the parameter $a$.

\textbf{Definition 1.1}
$$
T(bD_q)\{f(a)\}=\sum_{n=0}^{\infty}\frac{b^{n}}{(q;q)_n}\{{D_q}^n\{f(a)\}\},
$$
and
$$
 E(b\theta)\{f(a)\}=\sum_{n=0}^{\infty}\frac{q^{n(n-1)/
2}b^{n}}{(q;q)_n}\{{\theta}^n\{f(a)\}\}.
$$

The first and the second Cauchy operators and the operators
$T(b\theta)$ and $E(bD_q)$ below also act in this way.

When using the methods of operators,  rearrangememts of series are
often emploied. We know that rearrangememts of series must be under
the condition that the series are absolutely convergent, which may
be not very easy to check sometimes. In his paper \cite{liu9}~~~~~,
Liu obtained the following two theorems using the $q$-functional
equation method.
\begin{thm}\label{1.1}
Let $f(a,b)$ be a two variables analytic function in a neighborhood
of $(a,b)=(0,0)\in \mathbb{C}^2 $, satisfying the $q$-difference
equation
$$
bf(aq,b)-af(a,bq)=(b-a)f(a,b).
$$
Then we have
$$
f(a,b)=T(bD_q)\{f(a,0)\}.
$$
\end{thm}
 \begin{thm}\label{1.2}Let $f(a,b)$ be a two variables analytic
function in a neighborhood of $(a,b)=(0,0)\in \mathbb{C}^2 $,
satisfying the $q$-difference equation
$$af(aq,b)-bf(a,bq)=(a-b)f(aq,bq).$$
Then we have
$$f(a,b)=E(b\theta)\{f(a,0)\}.$$\end{thm}

These two theorems give us a method to use operator $T(bD_q)$ and
$E(b\theta)$ without having to check the absolute convergence of the
series.

Originated by Liu's work, we give the above
 two theorems more general forms in the following section.

\section{ Four $q$-functional equations and the solutions}

\begin{thm}\label{2.1}
Let $f(a,b,c)$ be a three variables analytic
function in a neighborhood of $(a,b,c)=(0,0,0)\in \mathbb{C}^3 $,
satisfying the $q$-functional equation
\begin{equation}
c(f(a,b,c)-f(a,bq,c))=b(f(a,b,c)-f(a,b,cq)-af(a,bq,c)+af(a,bq,cq)).
\end{equation}
Then we have
$$f(a,b,c)=T(a,b;D_{q})\{f(a,0,c)\},$$
where $T(a,b;D_{q})$ is applied with respect to the parameter $c$.
\end{thm}

\begin{thm}\label{2.2}
 Let $f(a,b,c)$ be a three variables analytic
function in a neighborhood of $(a,b,c)=(0,0,0)\in \mathbb{C}^3 $,
satisfying the $q$-functional equation
\begin{equation}
c(f(a,bq,cq)-f(a,b,cq))=b(f(a,bq,cq)-f(a,bq,c)-af(a,b,c)+af(a,b,cq)).
\end{equation}
Then we have
$$ f(a,b,c)=T(-\frac{1}{a},ab;\theta)\{f(a,0,c)\},$$
where $T(-\frac{1}{a},ab;\theta)$ is also applied with respect to
the parameter $c$.
\end{thm}

 When letting $a\rightarrow 0$ in theorem 2.1 and
2.2, we get theorem 1.1 and 1.2, respectively. The proof of theorem
2.2 is similar to that of theorem 2.1 and so is omitted.

{\bf  Proof of Theorem 2.1 } We now begin to solve this equation.
From the theory of several complex variables, we assume that
\begin{equation}
f(a,b,c)=\sum_{n=0}^{+\infty}A_n(a,c)b^n.
\end{equation}
We substitute the above equation into (1) to get
\begin{eqnarray*}c\sum_{n=0}^{+\infty}(1-q^n)A_n(a,c)b^n&=&
(\sum_{n=0}^{+\infty}A_{n}(a,c)-\sum_{n=0}^{+\infty}A_{n}(a,cq)
 -a\sum_{n=0}^{+\infty}A_{n}(a,c)q^{n}
 \\&&+a\sum_{n=0}^{+\infty}A_{n}(a,cq)q^{n})b^{n+1}.
      \end{eqnarray*}
This is
\begin{eqnarray*}
c\sum_{n=1}^{+\infty}(1-q^n)A_n(a,c)b^n=
\sum_{n=1}^{+\infty}(1-aq^{n-1})(A_{n-1}(a,c)-A_{n-1}(a,cq))b^{n}.
      \end{eqnarray*}
Comparing the coefficients of $b^n$ gives
$$
c(1-q^n)A_n(a,c)=(1-aq^{n-1})(A_{n-1}(a,c)-A_{n-1}(a,cq)).
$$
Then we have
\begin{eqnarray*}
A_n(a,c)&=&\frac{(1-aq^{n-1})}{(1-q^n)}
\frac{(A_{n-1}(a,c)-A_{n-1}(a,cq))}{c}
\\&=&\frac{(1-aq^{n-1})}{(1-q^n)}D_q\{A_{n-1}(a,c)\},
\end{eqnarray*}
where $D_q$ is applied with respect to the parameter $c$.
Iterate
the above equation to get
\begin{eqnarray}
A_n(a,c)=\frac{(a;q)_n}{(q,q)_n}
D_q\{A_{0}(a,c)\}.
\end{eqnarray}
It remains to calculate
$A_{0}(a,c)$.
Putting
$b=0$
in (3), we immediately deduce that
$A_{0}(a,c)=f(a,0,c)$.
Substituting (4) back into (3)~gives
$$
f(a,b,c)=\sum_{n=0}^{\infty}
\frac{(a, q)_n}{(q;q)_n}(bD_q)^{n}\{f(a,0,c)\}=T(a,b;D_{q})\{f(a,0,c)\}.
$$
This completes the proof.

Using the $q$-functional equation method we easily obtain the
following two theorems.
\begin{thm}
\label{2.3} Let $f(a,b)$ be a two variables analytic
function in a neighborhood of
$(a,b)=(0,0)\in \mathbb{C}^2 $,
satisfying the $q$-functional equation
$$
af(a,b)+bf(aq,bq)=(a+b)f(a,bq).
$$
Then we have
$$
f(a,b)=E(bD_q)\{f(a,0)\},
$$
where the operator $E(bD_q)$ is defined by
$$
 E(bD_q)=\sum_{n=0}^{\infty}\frac{q^{n(n-1)/
2}(bD_q)^{n}}{(q;q)_n}.$$
\end{thm}
\begin{thm}\label{1.1} Let $f(a,b)$ be a two variables analytic
function in a neighborhood of $(a,b)=(0,0)\in \mathbb{C}^2 $,
satisfying the $q$-functional equation
$$
af(aq,bq)+bf(a,b)=(a+b)f(aq,b).
$$
Then we have
$$f(a,b)=T(-b\theta)\{f(a,0)\},$$
where the operator $T(b\theta)$ is defined by
$$
T(b\theta)=\sum_{n=0}^{\infty}\frac{(b\theta)^{n}}{(q;q)_n}.
$$
\end{thm}

From all the theorems in this little paper,
we can look at the q-operators from a different
standpoint.

\end{document}